\documentclass[a4paper,12pt]{article}

\usepackage{citesort,geometry,enumerate,amsmath,amssymb,latexsym}                     
\def\baselinestretch{1.1}
\newtheorem{example}{Example}[section]
\newtheorem{defn}[example]{Definition}
\newtheorem{prop}[example]{Proposition}
\newtheorem{thm}[example]{Theorem}

\newenvironment{pf}{{\bf Proof:}}{$\Box$

\def\leq{\leqslant}

\mbox{}}
\begin{document}
\title{\Large \bf   SOME RESULTS FOR THE  LOCAL SUBGROUPOIDS 
 }

\small{ \author  {
\.{I}lhan \.{I}\c{c}en    \\  University of  \.{I}n\"{o}n\"{u} \\
Faculty of Science and Art
\\ Department of
Mathematics
\\ Malatya/ Turkey \\ iicen@inonu.edu.tr \\ \and  Osman Mucuk 
\\  University of  Erciyes \\
Faculty of Science and Art
\\ Department of
Mathematics
\\ Kayseri/ Turkey \\ mucuk@kayseri.edu.tr}}

\maketitle

\begin{abstract}
The notion of local subgroupoids as generalition of a local 
equivalence relations was defined  by the first author and R.Brown. 
Here we investigate some relations between transitive components and
coherence  properties of the local subgroupoids.

\end{abstract}

\section*{Introduction}
	Local to global problems play a very important role in mathematics. The most important concept 
in this context is a sheaf on a topological space. A sheaf is a way 
of describing a class of functions, sets, groups, etc
\cite{Sw,Ma-Mo,Bre}. For instance, a class 
of continuous functions on a topological space $X$ is very important in sheaf theory. 
The description tells the way in which a function $f$ defined on an open 
subset $U$ of $X$ can be restricted to functions $f|V$ on open subsets 
$V\subseteq U$ and then can be recovered by piecing together the restrictions
to the open subsets \cite{II1}. This applies not just to functions, but also to other mathematical 
structures  defined `locally' on a space $X$. Another set of local 
descriptions has been given in the notion of foliation (due to 
Ehresmann ) and also in the notion of local equivalence relation (due to 
Grothendieck and Verdier). 

The concept of the local equivalence relations,  
which was  introduced by Grothendieck and Verdier \cite{Gr-Ve} 
in a series of exercises presented as open problems 
concerning the construction of a certain kind of topos 
was investigated further by Rosenthal  \cite{Ro1,Ro2} 
and more recently by  
Kock and Moerdijk  \cite{Ko-Mo1,Ko-Mo2}. 
A local equivalence relation is a global section of a sheaf ${\cal E}$ 
defined by the presheaf  $E = \{ E(U), E_{UV}, X \}$, 
where $E(U)$ is the set of all equivalence relations on 
the open subset $U$ of $X$
and $E_{UV}$ is the restriction map from $E(U)$ to $E(V)$, for $V\subseteq U$. Moreover this presheaf
is not a sheaf. The key idea in this case is connectedness of the equivalence 
classes.

	An equivalence  relation on a set $U$ is just a wide  
subgroupoid of the indiscrete groupoid $U\times U$ on $U$. 
Thus it is  natural to consider 
the generalition which replaces the indiscrete groupoid on the 
topological space on $X$ by any groupoid $G$ on $X$. So we define 
a {\it local subgroupoid} of a groupoid $G$ on a topological space $X$
as a global section of the sheaf ${\cal L}$ associated to the presheaf
$L_G = \{  L(U),  L_{UV}, X \}$,
where $L(U)$ is the set  of all wide subgroupoids of $G|U$ 
and $L_{UV}$ is the restriction map from $L(U)$ to $L(V)$
for $V\subseteq U$. The concept of  local subgroupoids and their 
properties  were given  in \cite{Br-II1,BIM,II} .

In this paper, we  obtain the relation between the local subgroupoids 
and the local equivalence relations. 
Many coherence properties of local equivalence relations are 
investigated for local subgroupoids and  
obtained a topological foliation from the local subgroupoids. 

\noindent
{\bf Acknowledgment:} We would like to thank Prof. Ronald Brown 
for introducing us  to this area  and his help and encouragment.

\section{Local Subgroupoids}
Suppose $X$  is   a topological space   and   $G$ is a  groupoid on $X$.
For any open subset $U$  of $ X$, we write  $L_G(U)$  for 
the set of all wide subgroupoids of $ G|U$, 
where $G|U$ is the full subgroupoid of $ G$ on $U$.
Let $V, U\subseteq X$ be open sets such that $V\subseteq U$. 
If $H$ is a wide subgroupoid of $G|U$, then $H|V$ is a wide subgroupoid 
of $G|V$.
So there are restriction maps 
\[  L_{UV}\colon L_G(U) \longrightarrow L_G(V), \ \ H \longmapsto H|V \]
which define a presheaf $L_G \colon {\cal O}(X)^{op}\longrightarrow
Sets.$

However $ L_G$ is not a sheaf  in general \cite{Ro1,BIM}. 
It is well-known that the presheaf $L_G$ defines a sheaf  
$p\colon {\cal L}_G\rightarrow X$ \cite{Bre,II,Ma-Mo,Sw}.
For an open set $U$ in $X$, let $\Gamma (U, {\cal L}_G)$  be 
the set of sections. 
Then every element $H\in L_G(U)$ is associated with a 
section $s\in \Gamma(U,{\cal L}_G)$. If 
$x\in X$ and $\sigma \in p^{-1}(x)={\cal L}_x$, 
then there is an open 
neighbourhood  $U$ of $x$ in $X$ and an $s \in \Gamma (U,{\cal L}_G)$ such that
\[      \sigma = (U, H)_x = s(x).     \]
Such an element is called a {\it germ} at $x$. Similarly 
one can define a global section $s:X\to {\cal L}_G$  of the sheaf
$p\colon {\cal L}_G\rightarrow X$. 
\begin{defn}{\rm    
A {\it local subgroupoid} of $G$ on the topological space $X$ is a global  section  of  the 
sheaf ${\cal L}_G$ associated to the presheaf $L_G$.}
\end{defn}

Note that a local subgroupoid $s$ of $ G$  may be  given 
by local date called {\it atlas} \cite{II}; An atlas 
${\cal U}_s = \{(U_i, H_i): i\in I\}$ for the  local subgroupoid $s$ of
$G$ consist of an open cover  ${\cal U} =\{U_i:i\in I\}$ of $X$ and 
for each $i\in I$ a wide subgroupoid $H_i$  of $G|{U_i}$ such
that for all $x\in X, i\in I,$ if $x\in U_i$, then $s(x)=(U_i,
H_i)_x$.

If $s$ is a local subgroupoid of $G$ defined by an atlas 
${\cal U}_s = \{(U_i, H_i): i\in I\}$ and $V$ is an open subset of $X$ 
then $s|V$ is the local subgroupoid of $G|V$ defined by the atlas 
${\cal U}_s\cap V = \{(U_i\cap V, H_i|(U_i\cap V)): i\in I\}$. It is easy to
verify this is an atlas and as a section $s|V$ is just the 
restriction of $s$  to the open subset $V$.

Now  we give some examples of local subgroupoids
( see also \cite{BIM,II}).

\subsection{Examples}

\begin{example}{\rm 
Let $X$ be a topological space.
Every local equivalence relation on $X$ is a local subgroupoid  
of $X\times X$.
 For an  open set $U$ in $X$, let $E(U)$ be the set of  all 
 equivalence
relations on $U$. 
This gives us a presheaf $E : {\cal O}(X)^{op} \rightarrow Sets $  and so
it defines a sheaf ${\cal E}$ on $X$.
We obtain a local equivalence relation $r$ of global section of 
${\cal E}$.Obviously $E(U)$ is the set of all wide 
subgroupoid of $U\times U$. Hence the local equivalence relation 
$r$ on $X$ is a local subgroupoid of
$X\times X$. }
\end{example}
\begin{example}
{\em Any topological space $X$ can be considered as a 
groupoid on itself in which every element is  identity \cite{Br,Ma}.
We define a local subgroupoid $s$ of $X$  by an atlas
 ${\cal U}_s=\{(U_x, U_x) : x\in U_x \}$ where $U_x$ is open in $X$.}
\end{example} 
The following two examples give us the relations between local subgroupoids 
and local equivalence relations.
\begin{example}{\em 
Let 
${\cal X} = \{ X_i \mid i\in I \}$ be a
partition of a topological space $X$  and let $R$ be an equivalence 
relation defined by the partition. Then for a group $K$ , $G = R\times K$ becomes a groupoid on $X$
with the composition 
$((x, y), k_1)((y, z), k_2) =  ((x, z), k_1k_2))$
such that source and target maps $\alpha, \beta$ are  
the canonical projections 
on $R$ 

For open set $U_i$ of  $X$, 
let $E(U_i)$ be the set of  equivalence relations on $U_i$
and similarly
$L(U_i)$ the set of wide subgroupoid of $G|{U_i}$.
Then  $E$ and  
$L$ define presheaves on $X$. We obtain a natural transformation
$\gamma \colon E\rightarrow L$,   
by $\gamma_U\colon E(U_i)\rightarrow L(U_i)$, $R_i\mapsto R_i\times K$,
for each open set $U_i\subseteq X$, which 
defines  a presheaf morphism  $\gamma \colon E\rightarrow L$ and so 
it gives 
rise to a sheaf morphism  $\gamma^* \colon {\cal E}_R \rightarrow {\cal L}_G$.

If $r$ is a local equivalence relation of $R$ given by an atlas 
${\cal U}_r = \{ (U_i, R_i) \mid i\in I\}$, then $\gamma^* r$ becomes  a
local subgroupoid of $G$ on $X$ defined by  an atlas 
${\cal U}_{\gamma^* r} = \{ (U_i, R_i\times K) \mid i\in I \}$. In other word, a 
local equivalence relation of an equivalence relation $R$ on a topological 
space $X$ defines a local subgroupoid of the groupoid  $R\times K$ on $X$, 
while $K$ is a group.}
\end{example}
\begin{example} {\em 
Let $G$ be a groupoid on a topological space $X$. 
We  obtain a local
equivalence relation $r$ on $X$ to be a global section  of the sheaf ${\cal E}_{X\times X}$
associated with the presheaf $E_{X\times X}$
and a local subgroupoid $s$ of $G$ is a global section of the sheaf ${\cal L}_G$
associated with $L_G$. The anchor map 
\[ \Upsilon = [\alpha, \beta] \colon G \longrightarrow X\times X, 
\ \ \ g\mapsto (\alpha(g), \beta(g))   \]
 defines a  morphism of presheaves $\Upsilon \colon L(U)\rightarrow E(U)$, $H\mapsto R$. 
 Thena natural transformation
$\Upsilon \colon L_G\longrightarrow E_{X\times X}$ 
define a sheaf morphism
$\Upsilon^*\colon{\cal L}_G\longrightarrow {\cal E}_{X\times X}$
 such that $\Upsilon^*(s)=r$.}
\end{example}

\section{Coherent Local Subgroupoids}

We first recall  some elementary but essential basic 
facts from \cite{Br-II1,BIM,II}.

The set $L_G(X)$ of all wide subgroupoids of $G$ is a poset under
inclusion. We write $\leq$ for this partial order.

Let ${\bf Loc}(G)$ be the set of local subgroupoids of $G$. We define
a partial order $\leq$ on ${\bf Loc}(G)$ as follows.

Let $x\in X$. We define a partial order on the stalks
${p}^{-1}(x)={\cal L}_x$ by $(U', H')_x\leq (U, H)_x$ if there is
an  open neighbourhood $W$ of $x$ such that $W\subseteq U\cap U'$
and $H'|W$ is a subgroupoid of $H|W$. Clearly this partial order
is well defined. It induces a partial order on ${\bf Loc}(G)$ by
$s\leq t$ if and only if $s(x)\leq t(x)$ for all $x\in X$.

We now fix a groupoid  $G$  on $X$, so that $L_G(X)$ is the set of wide
subgroupoids of $G$, with its inclusion partial order, which we
shall write $\leq$.

We define poset morphisms
\[ loc_G : L_G(X)\to {\bf Loc}(G)    \ \ \mbox{and}\
\ \ \ \  glob_G:{\bf Loc}(G)\to L_G(X) \]
as follows.
We abbreviate $loc_G$, $glob_G$ to $loc, glob$.
\begin{defn}{\rm
If $H$ is a wide subgroupoid of the groupoid $G$ on $X$, then  $loc(H)$ is the
local subgroupoid defined by
\[             loc(H)(x) =(X, H)_x.    \]
Let $U$ be an open subset of $X$. Then we have notions of local
subgroupoids of $G|U$ and also of the restriction $s|U$ of a local
subgroupoid $s$ of $G$. Clearly if $H$ is a wide subgroupoid of
$G$ then $loc(H|U)= (loc(H))|U$.

Let $s$ be a local subgroupoid of  $G$. Then $glob(s)$ is the wide
subgroupoid of $G$ which is the intersection  of all wide
subgroupoids $H$ of $G$ such that $s\leq loc(H)$.} \end{defn}

We think of $glob(s)$ as an approximation to $s$ by a global
subgroupoid.

For any wide subgroupoid $H$ of $G$, $glob(loc(H))\leq H$.
However, $s\leq loc(glob(s))$ need not hold. Examples of this are
given in Rosenthal's paper \cite{Ro1} for the case of local
equivalence relations.
 We therefore adapt from \cite{Ro1,Ro2} some
notions of coherence.

\begin{defn}
Let $s$ be a local subgroupoid of $G$ on $X$.
\begin{enumerate}[(i)]
  \item   $s$ is called {\it coherent} if  $s\leq loc(glob(s))$.
\item  $s$ is called {\it totally coherent} if for every open set $U$ of $X$, $s|U$
is coherent.
\item $s$ is called {\it globally coherent} if $s = loc(glob(s))$.
\end{enumerate}

\end{defn}
\begin{defn}
Let $H\in L_G(X)$, so that $H$ is  a wide subgroupoid of $G$.
\begin{enumerate}[(i)]
  \item  $H$ is called {\it locally coherent} if $loc(H)$ is coherent.

  \item  $H$ is called {\it coherent} if $H = glob(loc(H))$.
\end{enumerate}
\end{defn}

Coherence of $s$ says that in passing between local and global
information nothing is lost due to collapsing. Notice also that
these definitions depend on the groupoid $G$.

The next proposition gives an alternative description of $glob$.

Let ${\cal U}_s = \{(U_i, H_i):i\in I\}$ be an atlas for the local
subgroupoid $s$. Then $glob({\cal U}_s)=H_{\cal U}$ is defined to be the
subgroupoid of $G$ generated by all the $H_i, i \in I$.

An atlas ${\cal V}_{s}=\{(V_j, {s}_j): j\in J\}$ for $s$ is
said to refine ${\cal U}_s$ if for each index $j\in J$
there exists an index $i(j)\in I$ such that $V_j\subseteq U_{i(j)}$
and $s_{i(j)}|V_j = s_j$ and written as ${\cal V}_{s}\leq {\cal U}_s$.

\begin{prop} \label{refin2}
Let  $s$  be a local subgroupoid of $G$ given by the atlas ${\cal
U}_s = \{(U_i, H_i):i\in I\}$. Then $glob(s)$ is the intersection
of the subgroupoids $glob({\cal V}_s)=H_{\cal V}$ of $G$ for all refinements
${\cal V}_s$ of ${\cal U}_s$, i.e., 
\[  glob(s) = \cap \{H_{\cal V} : {\cal V}_s\leq{\cal U}_s \}. \]
\end{prop}
\begin{pf}
Let $K$ be the intersection given in the proposition.

Let $H$ be a subgroupoid of $G$ on $X$ such that $s\leq loc(H)$.
Then for all $x \in X$ there is a neighbourhood $V$ of $x$ and
$i_x \in I$ such that $x \in U_{i_x}$ and $ H_{i_x}|V_x \cap
U_{i_x} \leq H$. Then ${\cal W}= \{(V_x \cap U_{i_x},H_{i_x}|V_x
\cap U_{i_x}) : x \in X\} $ refines ${\cal U}_s$ and $glob({\cal
W}) \leq H$. Hence $K \leq H$, and so $ K \leq glob(s) $.

Conversely, let ${\cal V}_s = \{ (V_j,H'_j) : j \in J \}$ be an
atlas for $s$ which refines ${\cal U}_s$. Then for each $ j \in J$
there is an $i(j) \in I$ such that $ V_j \subseteq
U_{i(j)},H'_j=H_{i(j)}|V_j$. Then $s \leq loc(glob({\cal V}_s))$. Hence
$glob(s) \leq glob({\cal V}_s)$ and so $glob(s) \leq K$.
\end{pf}

It is easy to show that if
for every open cover ${\cal V} $ of $X$,  $H\in L_G(X)$   is
generated by the subgroupoids $ H|V , V \in {\cal V} $ then 
$H$ is coherent. In fact,  atlas ${\cal U}_{loc(H)}= \{(X, H)\}$ 
which is refined by ${\cal V}_H= \{(V,H|V): V \in \cal V\}$  
for any open cover $\cal V$, hence 
\[  glob(loc(H)) = \cap \{H_{\cal V} : {\cal V}_s\leq{\cal U}_s \}=H. \] 

\begin{example}{\rm 
If $X$ is a topological space then its fundamental groupoid
$\pi_1X$ is coherent. Let ${\cal V}$ be an open cover of $X$. 
Let $[a]\in \pi_1X$, where
$a$ is a path  in $X$. Then by Lebesgue Covering Lemma, 
$a= a_1+\cdots +a_n$ where $a_i$ is a
path in some  $V_i\in {\cal V}$. So $\pi_1X$ is generated by
$(\pi_1X)|V$, $V\in {\cal V}$. }
\end{example}
\begin{example}{\em As we mentioned earlier, any topological space 
$X$ can be considered as a groupoid on itself and its local
subgroupoid is given an atlas $\{(U_x, U_x): x\in U_x\}$. 
Clearly $loc(X) = s$, i.e., 
$loc(X) (U) = (U, X\mid_U)=(U,U)$.
Moreover $glob(s) = X$, since $glob (s) = glop (loc (X)) = X$.
So $X$ is coherent and $s$ is globally coherent.}
\end{example}

\begin{example}{\em 
A bundle of groups $p\colon G\rightarrow X$ is a groupoid  whose source and target 
maps are equal, i.e. $\alpha = \beta = p$. 
For  an open set $U$ in $X$, $L_G(U)=\{G|U\}$ and
for open sets $V,U$ in $X$ with $V\subseteq U$,
restriction morphism 
$L_{UV}\colon L_G(U) \longrightarrow L_G(V), \ \  G_U \mapsto G_U\mid_V
$ define the presheaf $L_G$ 
and so a sheaf ${\cal L}_G$ as usual way.
So let $s$ be a local subgroupoid of the bundle of groups $G$ on $X$. 
Then $s$ is a globally coherent local
subgroupoid of  $G$ and $G$ is a coherent groupoid on $X$. In fact,
now, let ${\cal V} = \{ V_x  : x\in X\}$ be an open cover
of $X$ such that for each $x\in X$, $x\in V_x\subseteq U_x$, where
$ {\cal U} = \{ U_x : x\in X \}$ is also open cover  of $X$. Let $H_{\cal V}$ be the subgroupoid of $G$ generated
by $\{G\mid_{V_x}  : x\in X\}$.  $p^{-1} (U_x) = G\mid_{U_x}\supseteq H_x\mid_{V_x}$.
But $H_{\cal V} = G$ and  $glob(s) = \cap \{ H_{\cal V} : {\cal V}\leq
{\cal U} \} = H_{\cal V} = G$.
So $loc(glob(s)) = loc (G) = s$. Hence $s$ is a globally coherent. Since
$loc(G) = s$ and $glob(loc(G)) = glob(s) = G$, $G$ is coherent.}
\end{example}
More examples of the coherence
properties for local subgroupoids are given in  \cite{BIM}.

\subsection{Transitive components and coherence properties}

We give some results  which explain the relations between the 
transitive
components of $G$ and its local subgroupoids .
Recall that a groupoid $G$ on $X$ is called {\it transitive} if 
the set of morphism $G(x, y)$ from $x$ to $y$ is non-empty, 
for all $x, y\in X$. This defines an 
equivalence relation 
on $X$. The equivalence class containing $x\in X$ is denoted by $M_x$ 
and called {\it transitive components} of $G$ containing 
$x$ \cite{Br,Ma}.
\begin{prop}\label{Pr39}
 Let $s$ be a local subgroupoid of $G$ on 
$X$ such that  $ s = loc(H)$, for $H\in L_G(X)$. Then the 
transitivity componentsof $H_{\cal V}$ are relatively closed 
and open in the transitivity components of $H$.
\end{prop}
\begin{pf}
Let $M_{x,v}$ and $M_x$ denote the transitivity components of $x$ 
in $H_{\cal V}$ 
and $H$ respectively. Clearly,  $M_{x,v}\subseteq M_x$, because, 
$s = loc(H)$, $glob(s) = glob (loc(H))\subseteq H$, so 
$H_{\cal V}\subseteq H$.
Let $y \in M_{x,v}$. Then there is $g\in H_{\cal V}(y,x)$ with 
$g = h_n...h_2.h_1$
such that  $h_1\in H\mid_{V_1}(y,x_1),  h_2\in H|\mid_{V_2}(x_1, x_2),..., 
 h_n\in H\mid_{V_{n+1}}(x_n, x)$.
Take $V_1\cap M_x$ and let $z\in V_1\cap M_x$. Hence $h\in H(z,x)$ and since
$k\in H(y,x)$, $h^{-1} k \in H_{\cal V}(y, z)$. Thus, 
$h^{-1} k\in H_{\cal V}$ and $z\in M_{x,v}$.
Thus $z\in V_1\cap M_x\subseteq M_{x,v}$ and so $M_{x, v}$ is open in
$M_x$.
Now we have to  show that $M_{x,v}$ is closed in $M$. Let $z\in\overline{M_{x,v}}$
be the closure relative to $M_x$. For every open neighbourhood $U$ of $x$, we have 
$V_z$,  take an element $y\in V_z\cap M_{x,v}$. Then, there is a $g\in
H_{\cal V}(y,x)$ and since 
$h\in H(z, x)$, we have $h^{-1} g\in H(y,z)$. Since $y,z\in V_z$, 
 $h,g\in H_{\cal V}(y,z)$ and $z\in M_{x,v} $. 
 Thus $M_{x,v} = \overline{M_{x,v}}$
\end{pf}
\begin{thm}\label{Th310}
Let $s$ be a local subgroupoid of $G$ on 
$X$ such that  $ s = loc(H)$, for $H\in L_G(X)$.
Suppose that  for every $x\in X$
there is an open neighbourhood $W_x$ of $x$ such that $H\mid_{W_x}$ has connected transitivity  
components. Then $s$ is coherent.
\end{thm}
\begin{pf}
Suppose that $s$ is not coherent. Then, for some $a\in X$, 
we have 
$s(a)\not\leq loc(glob(s))(a)$, i.e., given any open neighbourhood 
$W$ of $a$,there is a cover $\{ V_x :x\in X\}$  and 
$y, z\in W$ such that there
exists an $h\in H(y,z)$, but  $h\not \in H_{\cal V}(y,z)$. 
In particular, this is true for $W_a$.
By Proposition \ref{Pr39}, the transitivity component of 
$y$ in $H\mid_{W_a}$ is open and closed
in the transitivity component of $y$ in $H|{W_a}$, 
which is connected.
This is a contradiction. Hence  $h\in H_{\cal V}(y,z)$.
\end{pf}

\begin{thm}\label{Th312}
Let $H\in L_G(X)$ and  $H$ has connected transitivity components.
Then $H = glob(loc(H))$. Conversely, if $H = glob(loc(H))$ and $H$ has closed transitivity 
components, then it has connected transitivity components.
\end{thm}
\begin{pf}
Given an open cover ${\cal V}=\{ V_x :x\in X\}$ of $X$. 
The groupoid $H_{\cal V}$
generated by  $\{H\mid_{V_x} \}$ is contained in $H$ 
and  by Proposition \ref{Pr39}, since transitivity components  of  
$H_{\cal V}$ are relatively open and closed  in those of $H$, 
which are connected, 
so we  must have $H_{\cal V} = H$, and hence  
$H = glob(loc(H)) = \cap\{ H_{\cal V}:{\cal V}\leq {\cal U}\}$.

If $H = glob(loc(H))$, for every cover ${\cal V}$, 
$H_{\cal V} = H$ .
Let $a\in X$  such that $M_a$, the transitivity component of $a$ in $H$, is 
not connected. Let $U$ and $V$  be open sets separating $M_a$.
Let ${\cal U} = \{ U, V, X-\{x\} \}$. Choose $x,y \in M_a$ such that 
$x\in U\cap M_a$,  $y\in V\cap M_a$. Then there exists 
$g\in H(x,y),$ but 
$g\notin H_{\cal V}(x,y)$ since $(U\cap M_a)\cup (V\cap M_a) = M_a$ and
they are disjoint, since $H_{\cal V}\not\subseteq H$ we have that  $glob(loc(H))\subseteq H$.
This is a contradiction.
\end{pf}
The following results  denote   relations between the previous
definitions.
\begin{prop}
i) Suppose that $s$ is globally and totally coherent of a groupoid $G$ on $X$. 
If $U$ is open 
in $X$, then $s|U$ is globally coherent.

ii) If there is an open cover ${\cal V}= \{ V_x : x\in X \}$ of $X$ such that
$s|{V_x}$ is globally and totally coherent for all $x\in X$, 
then $s$ is totally coherent.
\end{prop}
\begin{pf}
i) We  have  $s = loc(glob(s))$. By definition, 
$glob(s|U)\subseteq glob(s)|U$,
hence $loc(glob(s|U))\leq loc(glob(s)|U) = loc(glob(s))|U = s|U$.
Sine $s|U$ is coherent, by totally coherence of $s$, we have 
$s|U\leq loc(glob(s|U))$. So $s|U = loc(glob(s|U))$, i.e., 
$s|U$ is globally coherence.

ii) By (i), if $U$ is open in $X$ and $s|{V_x}$ is globally 
coherence for all 
$x\in X$, then $s|{U\cap V_x}$ is globally coherent. 
Thus  $s|{U\cap V_x}$ 
$= loc(glob(s))|{U\cap V_x} \leq loc(glob(s|U))|{V_x}$, since this holds for all $x\in X$,
we have $ s|U\subseteq loc(glob(s|U))$, i.e., $s$ is totally coherent.
\end{pf}

\subsection{Topological foliations}
One of Ehresmann's approaches to the foundations of foliation theory 
\cite{Eh} goes via the consideration of a topological space equipped
with a further `nice' topology. Such nice topologies appear also
in the context of local equivalence relations and have been
considered in \cite{AGV} and in \cite{Ro2}. We shall need the
following elaboration of this idea for the local subgroupoids.
Let $s$ be a local subgroupoid of $G$ on a topological space $X$ which 
is given by an atlas $\{(U_x, H_x): x\in X\}$.
We can define a new topology on $X$ denoted by $X^s$. 
The underlying set of $X^s$ is $X$. Let $M_{x,a}$ denote 
the transitivity 
components of $x$  in $U_a$ for the subgroupoid 
$H_a\in L_G(U_a)$. Let the topology 
of $X^s$  be generated by the $M_{x,a}$ , $x\in U_a$ and 
the open sets of $X$.
Then its basic open sets any set of the form $U\cap M_{x,a}$ 
where $U$ is open in $X$, thus this topology is the coarsest for
which the original open sets as well as  transitive component for
$H_a$ are open, and $X^s$ is topologically the disjoint union of the
transitive component for $H_a$, each of them with its subspace
topology from $X$.
Since the 
topology on $X^s$ is finer than that of $X$, 
$I\colon X^s\rightarrow X$, 
the identity map, 
is continuous. Hence $X^s$ is a {\it topological foliation} .
The notion of topological foliation was defined by Ehresmann
\cite{Eh}.
\begin{thm}\label{Th314}
Let $s$ be a coherent local subgroupoid of the groupoid $G$ on $X$. 
Then the transitivity components of $glob(s)$
are connected components of $X^s$.
\end{thm}
\begin{pf}
Let $H = glob(s)$. Since $s$ is coherent ($s\leq loc(H) $), for each $a\in X$,
choose an open neighbourhood $W_a\subseteq U_a$ such that 
$H_a|{W_a}\subseteq H|{W_a}$.
If $M$ is a transitivity component of $G$, we shall show that 
\[        M = \bigcup_{a\in M}(M_{a,a}\cap W_a).       \]
If $z\in M_{a,a}\cap W_a$ for some $a\in M$, then 
$h\in H_a(z,a)|{W_a}$ and hence $h\in G$.
Since $a\in M$ and $M$ is the transitivity component of $G$, 
then $z\in M$. Hence
$M$ is a union of open sets in $X^s$ and so is open.

We prove $M$ is closed in $X^s$, let $x\in {\overline M}$, closure is relative to $X^s$, 
 then $M_{x,x}\cap W_x$ meets $M$. Let $a\in (M_{x,x}\cap W_x)\cap M$. Then,
$k\in {H_x}|{W_x}(a,x)\subseteq G|{W_x}$. Since $a\in M, x\in M$.
Thus $M = \overline M$ and $M$ is closed in $X^s$.

Since $M$ is open and also closed, 
if it is transitivity connected, we have to show that
it is a connected component. Since $s$ is coherent and the 
topology of $X^s$ is finer than that of $X$, 
it follows that $s\leq loc(H)$ in $X^s$, from which it follows that
$glob(s)\leq glob(loc(H))$ and hence $glob(s) = glob(loc(glob(s))$, i.e.,
$H = glob(s)$ is coherent on $X^s$. Since its transitivity components are 
closed by Theorem \ref{Th312}, they are connected.
\end{pf}
As a result, some properties of local equivalence relations 
can be described by the results
centred around the notion of local subgroupoids. The interplay of the functors 
{\em glob} and {\em loc}  says a lot about  local subgroupoids on 
arbitrary topological spaces.

{}

\end{document}